\documentclass[reqno]{amsart}

\usepackage{amsthm,amssymb}

\usepackage{latexsym,abbreviations, indentfirst}
\usepackage{color}
\usepackage{epsfig}

\usepackage{natbib}

\bibpunct{(}{)}{;}{a}{,}{,}

\newtheorem{theorem}{Theorem}
\newtheorem{proposition}{Proposition}
\newtheorem{corollary}{Corollary}
\newtheorem{lemma}{Lemma}
\newtheorem{fact}{Fact}

\theoremstyle{definition}
\newtheorem{definition}{Definition}

\theoremstyle{remark}
\newtheorem{remark}{Remark}

\newcommand{\PPP}{\mathbb{P}}

\newcommand{\nid}{\noindent}

\renewcommand{\qedsymbol}{\ensuremath{\blacksquare}}

\begin{document}

\title{A note on recurrent random walks}
\author{Dimitrios Cheliotis}

\thanks{
Department of Mathematics, Bahen Center for Information
Technology, 40 St. George St., 6th floor, Toronto, ON, M5S 3G3,
Canada. \textit{E-mail address}: dimitris@math.toronto.edu,
\textit{URL}: http://www.math.toronto.edu/dimitris/}
\date{August 18, 2005}
\thanks{Research partially supported by an anonymous Stanford Graduate Fellowship
and a scholarship from the
 Alexander S. Onassis Public Benefit Foundation}

\maketitle

\begin{abstract}
 For any recurrent random walk $(S_n)_{n\ge1}$ on $\D{R}$, there
are increasing sequences $(g_n)_{n\ge1}$ converging to infinity
for which $(g_nS_n)_{n\ge1}$ has at least one finite accumulation
point. For one class of random walks, we give a criterion on
$(g_n)_{n\ge1}$ and the distribution of $S_1$ determining the set
of accumulation points for $(g_nS_n)_{n\ge1}$. This extends, with
a simpler proof, a result of K.L. Chung and P. Erd\"os. Finally,
for recurrent, symmetric random walks, we give a criterion
characterizing the increasing sequences $(g_n)_{n\ge1}$ of
positive numbers for which $\varliminf g_n|S_n|=0$.

\bigskip

\nid Keywords:random walk, recurrence, stable distributions,
symmetric distributions.
\end{abstract}

\section{Introduction} \label{RecNote}
We recall that a random walk in $\D{R}$ is any random sequence
$(S_n)_{n\ge1}$ with $S_n=X_1+\cdots+X_n$ and the $X_i$'s i.i.d.
having common distribution $F$. Such a walk is called recurrent if
for any $\gep>0$, with probability one, there are infinitely many
indices $n$ for which $S_n\in(-\gep, \gep)$. It is easy to see
that, for such a walk, it is possible to find a deterministic
sequence of positive numbers $(\gep_n)_{n\ge1}$ converging to zero
for which, with probability one, $S_n\in(-\gep_n, \gep_n)$
infinitely often. And in some cases it is possible to find a
deterministic sequence of positive numbers $(\gep_n)_{n\ge1}$
converging to zero for which, with probability one,
$S_n\notin(-\gep_n, \gep_n)$ eventually. This is the content of
the next proposition, whose proof is given in Section \ref{RecRW}.

\begin{proposition} \label{propRec}
Let $(Y_n)_{n\ge1}$ be a sequence of random variables, defined on
the same probability space, with $\varliminf_{n\to+\infty}|Y_n|=0$
a.s. Then

\begin{enumerate}
\item[(a)] There exists a sequence of positive numbers
$(g_n)_{n\ge1}$ increasing to infinity with
$\varliminf_{n\to+\infty}g_n|Y_n|=0$ a.s. \item[(b)] If
$\sum_{n=1}^{+\infty} \PPP(Y_n=0)<+\infty$, then there exists a
sequence of positive numbers $(g_n)_{n\ge1}$ increasing to
infinity with $\lim_{n\to+\infty}g_n|Y_n|=+\infty$ a.s.
\end{enumerate}
\end{proposition}

By determining the set of sequences $(g_n)_{n\ge1}$ with the
properties $g_n\nearrow +\infty$ and
$\varliminf_{n\to+\infty}g_n|S_n|=0$, we measure in some way the
``strength'' of the recurrence of the random walk. We would like
to establish a criterion characterizing these sequences in terms
of the distribution function $F$.

Our first theorem addresses this question for the case where $F$
is the distribution of a symmetric random variable. In such a
case, we call also $F$ symmetric.
\begin{theorem} \label{symmetric}
 Assume that $F$ is symmetric.
Then for any increasing sequence $(g_n)_{n\ge 1}$ of positive
numbers and $M>0$,
$$ \varliminf g_n|S_n| =\left\{
\begin{array}{ll} 0 \text{ a.s. }& \text{ if } \sum\limits_{n=1}^{+\infty}\PPP(g_n|S_n|<M)=+\infty,  \\
 +\infty  \text{ a.s. } & \text{otherwise.}
\end{array}
    \right.$$
\end{theorem}

For the statement of our next theorem, we define for any sequence
$x:=(x_n)_{n\ge1}$ of real numbers
$$K(x):=\{a\in\ol{\D{R}}:a \text{ is an accumulation point for }
(x_n)_{n\ge1}\},$$
 where $\ol{\D{R}}:=\D{R}\cup\{-\infty, +\infty\}$ is equipped with the usual topology.
We also recall the
following definition.
\begin{definition}\label{DomAttr}
The distribution $F$ belongs to the domain of attraction of a
distribution $R$ if there exist constants $a_n>0$, $b_n$ such that
the distribution of $(S_n-b_n)/a_n$ tends to $R$ as $n\to+\infty$.
\end{definition}

\begin{theorem} \label{stable}
Assume that $F$ belongs to the domain of attraction of a stable
law, with zero centering constants (i.e., in the definition above,
the $b_n$'s are zero), and moreover that the characteristic
function $ \varphi$ of $F$ satisfies Cramer's condition (C), i.e.,
$\varlimsup_{|t|\to+\infty}|\varphi(t)|<1$. Then for any
increasing sequence $(g_n)_{n\ge 1}$ of positive numbers, it holds
\begin{equation*}
K(gS)=
    \begin{cases} \{-\infty, +\infty\} & \text{a.s. if } \sum\limits_{n=1}^{+\infty}\frac{1}{a_n g_n}<+\infty,  \medskip \\

\ol{\D{R}}  & \text{a.s. if
}\sum\limits_{n=1}^{+\infty}\frac{1}{a_n g_n} =+\infty.
    \end{cases}
\end{equation*}
\end{theorem}

In the case where $F$ itself is the distribution of a symmetric
stable law with index in $[1,2]$, we have the following corollary.

\begin{corollary} For a symmetric stable law with index
$\ga\in[1,2]$ (i.e., characteristic function
$\varphi(t)=e^{-c|t|^\ga}$ for some $c>0$),
\[
K\big((n^\gep\,S_n)_{n\ge1}\big)=\left\{
\begin{array}{ll} \ol{\D{R}} &
    \text{ a.s. if } \gep\le 1-\frac{1}{\ga},  \medskip \\ \{-\infty, +\infty\}  & \text{ a.s. if }\gep>1-\frac{1}{\ga}. \end{array}
    \right.
\]
\end{corollary}

\begin{remark}
When $g_n=1$ for all $n$, Theorem \ref{symmetric} reduces to a
well known criterion for recurrence, which holds for all
distributions (see \cite{DU}, Chapter 3, Corollary 2.6)
\end{remark}

\begin{remark} If $(S_n)_{n\ge1}$
is a recurrent random walk and $(g_n)_{n\ge1}$ is an increasing
sequence of positive numbers, then $-\infty, +\infty\in K(gS)$
because $\varliminf_{n\to+\infty}S_n=-\infty,
\varlimsup_{n\to+\infty}S_n=+\infty$ (see \cite{DU}, Chapter 3,
Theorem 1.2). So the above theorems, in particular, characterize
the sequences $(g_n)_{n\ge1}$ for which $K(gS)=\{-\infty,
+\infty\}$, i.e., the ones that push the walk to infinity.
\end{remark}

\begin{remark}
Let $\varphi$ denote the characteristic function corresponding to
the distribution $F$. If $F$ is lattice with span $h>0$ and the
corresponding walk recurrent, then $\PPP(S_n=kh \text{ i.o.})=1$
for all $k\in\D{Z}$, and
$K(gS)=\{0\}\cup\{kh\lim_{n\to+\infty}g_n:k\in\D{Z}\setminus\{0\}\}$
for all increasing sequences $(g_n)_{n\ge1}$ of positive numbers.
In particular, $K(gS)\ne\{-\infty, +\infty\}$ for all such
sequences; i.e., the walk is very recurrent but for a trivial
reason. That is, because $\PPP(S_n=0 \text{ i.o.})=1$. The
remaining distributions, called nonlattice, are the ones
satisfying $|\varphi(t)|<1$ for all $t\in\D{R}\setminus\{0\}$. A
subclass of these are those that satisfy Cramer's condition (C),
and Theorem \ref{stable} is concerned with these. In that theorem,
the assumption that $F$ satisfies Cramer's condition (C) cannot be
weakened to $F$ nonlattice. To see this, take
$a\in\D{R}\setminus\D{Q}$ and the distribution $F$ that assigns
mass 1/4 to each of the numbers in the set $\{-1,1, -a, a\}$. $F$
belongs to the domain of attraction of the normal distribution,
with $a_n=\sqrt{n}$ and zero centering constants. Its
characteristic function is $\varphi(t)=(\cos t+\cos at)/2$ for all
$t\in\D{R}$ and satisfies $|\varphi(t)|<1$ for all
$t\in\D{R}\setminus\{0\}$ since $a$ is irrational. Thus, the
distribution is nonlattice. But it does not satisfy Cramer's
condition (C) because $\varlimsup_{k\to+\infty, k\in\D{N}}
(ak-[ak])=1$ due again to the irrationality of $a$; and
$\varphi(2k\pi)=(1+\cos((ak-[ak])2\pi))/2$ for $k\in\D{N}$. For
the sequence $(n)_{n\ge1}$ we have
$\sum\nolimits_{n=1}^{+\infty}(\sqrt{n}n)^{-1}<+\infty$. If our
theorem would apply to this case, we would have $\lim
n|S_n|=+\infty$ a.s. This is false, because for $n\ge1$ we have
$\PPP(S_{2n}=0)=4^{-2n}\sum\nolimits_{m=0}^n
(2n)!/[m!m!(n-m)!(n-m)!]\sim (\pi n)^{-1}$ (it is the same
calculation as for the two dimensional simple random walk since
$a$ is irrational), implying $\sum_{n=1}^{+\infty}
\PPP(S_{2n}=0)=+\infty$, and by a well known result,
$\PPP(S_{2n}=0 \text{ i.o.})=1$ (see \cite{DU}, Chapter 3, Theorem
2.2).
\end{remark}

The first result in the spirit of our theorems is the content of
Theorem 3 of \cite{CE}, where the authors characterize the
sequences $(g_n)_{n\ge 1}$  for which $(g_n/\sqrt{n})_{n\ge1}$
increases to infinity and $\varliminf_{n\to+\infty}g_n |S_n|=0$ .
They assume that the distribution of the $X_i$'s has finite
absolute fifth moment, zero mean value, and a non zero absolutely
continuous part. Later, V.V.Petrov (see \cite{PE}) improved their
result by assuming only finite second moment, zero expectation,
and Cramer's condition stated above; also he removed the
requirement that $(g_n/\sqrt{n})_{n\ge 1}$ converges to infinity
and assumed that it is just increasing.

Our second result is an extension of the work of Erdos, Chung,
Petrov as it determines all accumulation points of
$(g_nS_n)_{n\ge1}$ for a bigger class of distributions $F$ and
sequences $(g_n)_{n\ge1}$. The proofs of the aforementioned
authors are longer than ours because they establish from first
principles that a set of interest has probability one. In our
approach, we just prove that this set has positive probability,
and then we invoke the Hewitt-Savage zero-one law.

Finally, we should mention the following result of V. Petrov
(Theorem 6.21 in \cite{PE2}). Assume that the characteristic
function of the $X_i$'s satisfies Cramer's condition (C). Then for
any sequence $(g_n)_{n\ge1}$ of positive numbers,
$\sum_{n=1}^{+\infty} (\sqrt{n}g_n)^{-1}<+\infty$ implies
$\lim_{n\to+\infty} g_n|S_n|=+\infty$ a.s.

\section{Proof of the results} \label{RecRW}

In this section, we prove the two theorems and Proposition
\ref{propRec}. First, we give a version of the second
Borel-Cantelli lemma that we will use. It is a simple application
of the Kochen-Stone lemma (see \cite{DU}, Chapter 1, exercise
6.20), so we omit its proof.
\begin{lemma}\label{BC} Assume that $(A_n)_{n\ge 1}$ is a sequence
of measurable sets and that there are $n_0\ge 1$ and $c>0$ such
that $\PPP(A_j\cap A_k)\le c\,\PPP(A_j)\,\PPP(A_{k-j})$ for every
$k,j\ge n_0$ with $k-j\ge n_0$. Then
$$ \PPP(A_n \mbox{ i.o.})>0 \mbox{\quad iff \quad} \sum\limits_{n=0}^{+\infty}\PPP(A_n)=+\infty.$$
\end{lemma}
For $a\in\D{R}, \gep>0$,  and $n\ge 1$, we set
$A_n(a,\gep):=\left[g_n S_n \in (a-\gep, a+\gep) \right]$, and
observe that $$g_n S_n\in(a-\gep, a+\gep) \text{ i.o.  a.s. }
\Leftrightarrow \PPP(\varlimsup\limits_{n} A_n(a,\gep))=1.$$ Since
the $X_i$'s are i.i.d and $\varlimsup\limits_{n} A_n(a, \gep)$ is
an exchangeable event, by the Hewitt-Savage zero-one law, if it
has positive measure, then it has measure one. Both theorems
(\ref{symmetric} and \ref{stable}) are proved by showing that the
sequence ($A_n(a,\gep))_{n\ge1}$ satisfies the conditions of Lemma
\ref{BC} for all $\gep>0$ and appropriate $a$ (for $a=0$ for
Theorem \ref{symmetric}, and for all $a\in\D{R}$ for Theorem
\ref{stable}).

\nid To do this, we observe that for $a\in\D{R}, \gep>0$, and any
$j<k$,
\begin{multline*}
\PPP(A_j(a, \gep)\cap A_k(a, \gep))\\ =\PPP\Big(g_j
S_j\in(a-\gep, a+\gep), g_k(S_k-S_j)\in(-g_kS_j+a-\gep, -g_kS_j+a+\gep) \Big) \\
\le \PPP(A_j(a, \gep))\sup\limits_{y\in((a-\gep)/g_j,
(a+\gep)/g_j)}\PPP\Big( S_{k-j}\in(-y+\frac{a-\gep}{g_k},
-y+\frac{a+\gep}{g_k} )\Big).
\end{multline*}
And now we want to bound the last supremum by $c\,\PPP(A_{k-j}(a,
\gep))$ for some constant $c$ that may depend on $a, \gep$ but not
on $k, j$. It remains to establish such a bound under the
assumptions of either of the two theorems.

\medskip

\nid Theorem \ref{symmetric} is a consequence of the following two
results.

\begin{proposition} \label{symmetricprop}
 Assume that $F$ is symmetric.
Then for any increasing sequence $(g_n)_{n\ge 1}$ of positive
numbers and $M>0$,
$$ \PPP(g_n|S_n|< M \text{ i.o.})=\left\{
\begin{array}{ll} 1 & \text{ if } \sum\limits_{n=1}^{+\infty}\PPP(g_n|S_n|<M)=+\infty, \\
 0  & \text{ otherwise. }
\end{array}
    \right.$$
\end{proposition}

\begin{proof}
If $\sum\nolimits_{n=1}^{+\infty}\PPP(g_n|S_n|<M)<\infty$, we
apply the first Borel-Cantelli lemma. So assume

\nid $\sum\nolimits_{n=1}^{+\infty}\PPP(g_n|S_n|<M)=+\infty$. We
claim that
$\sum\nolimits_{n=1}^{+\infty}\PPP(g_{2n}|S_{2n}|<M)=+\infty$.
Indeed, for $n\ge1$ we have
\begin{multline*}\PPP(|S_{2n+1}|<M/g_{2n+1})=\int\PPP(|S_{2n}+y|<M/g_{2n+1})dF_X(y)\\ \le
4\PPP(|S_{2n}|<M/g_{2n+1})\le 4\PPP(|S_{2n}|<M/g_{2n}).
\end{multline*}
The first inequality follows from Lemma \ref{dominance} because
the symmetry of the distribution of the $X_i$'s implies that
$S_{2n}$ has a nonnegative characteristic function. The second
inequality follows from the monotonicity of $(g_n)_{n\ge1}$.

Combining
$\sum\nolimits_{n=1}^{+\infty}\PPP(g_{2n}|S_{2n}|<M)=+\infty$ with
the comments before the proof of this proposition (applied to the
random walk $(S_{2n})_{n\ge1}$ and the sequence
$(g_{2n})_{n\ge1}$) and Lemma \ref{dominance}, we get the desired
result.
\end{proof}

\begin{lemma} \label{dichotomy}
Assume that $F$ is symmetric and  $(g_n)_{n\ge 1}$ is an
increasing sequence of positive numbers. Then
$\sum\limits_{n=1}^{+\infty}\PPP(g_n|S_n|<M)$ is either infinite
for all $M>0$ or finite for all $M>0$.
\end{lemma}
\begin{proof}
Assume the contrary. Then there is an $a>0$ so that the above sum
is finite for $M<a$ and infinite for $M>a$. Take $\gep\in(0, a)$.
Since
$$\PPP(g_n|S_n|<a+\gep)=\PPP(g_n|S_n|\le a-\gep)+\PPP(g_n|S_n|\in (a-\gep,
a+\gep)),$$ it follows that \begin{equation}\label{infinsum}
\sum\limits_{n=1}^{+\infty}\PPP(g_n|S_n|\in(a-\gep,
a+\gep))=+\infty.\end{equation} As in the proof of Proposition
\ref{symmetricprop}, we will show that
$\sum\limits_{n=1}^{+\infty}\PPP(g_{2n}|S_{2n}|<\gep)=+\infty$.
Indeed, for $n\ge1$,
\begin{multline*} \PPP(g_{2n+1}|S_{2n+1}|\in(a-\gep, a+\gep)) \\ =\int
\PPP(S_{2n}\in(-x+\frac{a-\gep}{g_{2n+1}},
-x+\frac{a+\gep}{g_{2n+1}})\cup (-x-\frac{a+\gep}{g_{2n+1}},
-x-\frac{a-\gep}{g_{2n+1}}))\ dF(x)\\ \le
8\PPP(S_{2n}\in(-\frac{\gep}{g_{2n+1}}, \frac{\gep}{g_{2n+1}}))\le
8\PPP(S_{2n}\in(-\frac{\gep}{g_{2n}}, \frac{\gep}{g_{2n}})).
\end{multline*}
Also $\PPP(g_{2n}|S_{2n}|\in(a-\gep, a+\gep))\le 8
\PPP(g_{2n}|S_{2n}|<\gep)$. Thus, \eqref{infinsum} and the last
two inequalities imply
$\sum\limits_{n=1}^{+\infty}\PPP(g_{2n}|S_{2n}|<\gep)=+\infty$. A
contradiction since $\gep<a$.
\end{proof}

\medskip
The essential ingredient for the proof of Proposition
\ref{symmetricprop} and Lemma \ref{dichotomy} is the following
lemma.
\begin{lemma} \label{dominance}
If $X$ is a random variable with real and nonnegative
characteristic function $ \varphi$, then for all $x\in\D{R}$ and
$\gd>0$ we have
$$\PPP(X\in[x-\gd, x+\gd]) \le 4\PPP(X\in(-\gd, \gd)).$$
\end{lemma}

\begin{proof} Clearly, we can assume that $\gd=1$.
 Let $A:\D{R}\to\D{R}$ be defined by $A(x)=\int_0^1 \PPP(X\in[x-r,
x+r])\ dr$ for every $x\in\D{R}$. Then using the inversion formula
and the bounded convergence theorem, we get $A(x)=(2\pi)^{-1}
\int_{\D{R}}e^{-itx}(1-\cos t)t^{-2} \varphi(t) \ dt$, which
attains global maximum at $x=0$ because $\varphi$ is nonnegative.
Thus $\PPP(X\in(-1, 1))\ge A(0)\ge A(x)\ge \PPP(X\in[x-1/2,
x+1/2])/2$. Consequently $\PPP(X\in[x-1, x+1])\le 4
\PPP(X\in(-1,1))$.
\end{proof}

For the proof of Theorem \ref{stable}, we will use the following
two facts.
\begin{fact} Let $S_n, a_n, b_n$ be as in Definition \ref{DomAttr}, with
the law of $(S_n-b_n)/a_n$ converging to a stable distribution
with index $\ga$, not concentrated at zero. Then
$h(n):=a_n/n^{1/a}$ is a function slowly varying at infinity.
\end{fact}
\begin{fact}
If $h$ is a function slowly varying at infinity and $\gd>0$, then
$n^{-\gd}<h(n)<n^{\gd}$ for all big $n$.
\end{fact}
Fact 1 is contained in Theorem 2.1.1 of \cite{IB}. Fact 2 is Lemma
2 of \S 8, Chapter VIII in \cite{F}.

\medskip
\nid \textbf{Proof of Theorem \ref{stable}:}

Call $F_\ga$ the distribution function of the limiting stable law.
Its characteristic function is absolutely integrable (see
\cite{F}, Chapter XVII, Section 6), so it has a continuous bounded
density, call it $f_\ga$. We will consider two cases depending on
the value of $\ga$. The nontrivial is the second one.

\nid CASE 1: $0<\ga<1$.

\nid By Facts 1,2, we have $\sum_{n=1}^{+\infty}1/a_n <+\infty$,
so the series in the statement of the theorem converges always. It
is enough to prove that $S_n$ is transient. Let $M$ be the bound
in the density of $f_\ga$. By Lemma 2 in \cite{ST}, there exists
an $n_0\ge1$ and $h_0>0$ such that for all $x\in\D{R}$ and
$h\in[a_n^{-1}, h_0]$ we have
$$h(f_\ga(x)-1)\le \PPP(S_n/a_n\in [x, x+h])\le h(f_\ga(x)-1).$$
To apply that lemma one needs to have $F$ nonlattice, which is
true because $\varphi$ satisfies Cramer's condition (C). Now pick
$\gep\in(\ga,1)$. Then for $h:=2a_n^{-\gep}$ and $x:=-h/2$, we
have $\PPP(|S_n|<a_n^{1-\gep})=\PPP(|S_n/a_n|<a_n^{-\gep})\le
(M+1) a_n^{-\gep}$ for big $n$. Since $\gep/\ga>1$, Facts 1,2
imply that $\sum_{n=1}^{+\infty}
 a_n^{-\gep}<+\infty$, and the first Borel-Cantelli lemma gives
 that $|S_n|\ge a_n^{1-\gep}$ eventually. So the walk is transient since
 $1-\gep>0$.

\nid CASE 2: $1\le \ga\le 2$.

 \nid Assume that
$\varlimsup_{|t|\to+\infty}|\varphi(t)|=\theta<1$, and pick an
$\gep_0>0$ so that $\theta+\gep_0<1$.

\textbf{Step 1:}  First we prove the theorem for all sequences
$(g_n)_{n\ge1}$ as in the statement of the theorem that moreover
satisfy
\begin{equation}\label{CramerControl}
\lim_{n\to+\infty}g_n a_n(\theta+\gep_0)^n=0.
\end{equation}
Using the explicit expressions for $f_\ga$ (see Lemma 1 in
\cite{F}, Chapter XVII, Section 6), we can see that $f_\ga(0)>0$
for $\ga \in(1,2]$. For $\ga=1$, we use the fact that $F_\ga$ is
strictly stable (i.e. the centering constants are zero) to obtain
that its characteristic function is of the form
$\varphi_1(t)=e^{itc-d|t|}$ for all $t\in\D{R}$, where $c
\in\D{R}$ and $d>0$ are some constants (e.g., by exploiting the
relation $\varphi_1(nt)=(\varphi_1(t))^n$ for all $n\in\D{N}$ and
$t\in\D{R}$, and the general form of the characteristic function
of a stable distribution with index 1 given in Theorem 9.32 of
\cite{BREI}). Thus, $F_1$ is the distribution function of $Y d+c$,
where $Y$ has a Cauchy distribution, which implies that
$f_1(0)>0$.

By the Local Limit Theorem of Stone (see \cite{ST}), for fixed
$c>1$, there is an $h_0>0$ and $n_0>1$ such that for any interval
$I\subset[-h_0, h_0]$ and $n\ge n_0$ with $|I|>(\theta
+\gep_0)^n$, we have
\begin{equation}\label{StoneLocal}
\frac{1}{c}\,f_\ga(0)<\frac{1}{|I|}\,\PPP(\frac{S_n}{a_n} \in
I)<c\,f_\ga(0).
\end{equation}
This theorem applies because $\varphi$ satisfies Cramer's
condition (C).

Now take $a\in \D{R}$ and $\gep>0$. Since $g_n
a_n(\theta+\gep_0)^n\to 0$ and $a_n g_n \to +\infty$
 as $n\to+\infty$ (see Facts 1,2 regarding $a_n$),
 we can assume that for the above $n_0$ we also have $(|a|+\gep)/a_n g_n<h_0/2$,\, $g_n
 a_n(\theta+\gep_0)^n/2\gep<1$ for all $n\ge n_0$.
Thus, for $n\ge n_0$ and $I:=\big((a-\gep)/a_ng_n, (a+\gep)/a_n
g_n\big)$, \eqref{StoneLocal} applies and gives
 \begin{equation}\label{localEst}
 c^{-1}f_\ga(0)<\frac{a_n g_n}{2\gep}\,\PPP(A_n(a,\gep))<
 c\,f_\ga(0).
 \end{equation}
\nid Also
\begin{multline}
\PPP\Big( S_{k-j}\in(-y+\frac{a-\gep}{g_k}, -y+\frac{a+\gep}{g_k}
)\Big)\\ =\PPP\Big(
\frac{S_{k-j}}{a_{k-j}}\in(-\frac{y}{a_{k-j}}+\frac{a-\gep}{g_ka_{k-j}},
-\frac{y}{a_{k-j}}+\frac{a+\gep}{g_k a_{k-j}} )\Big)\\ \le
\PPP\Big(
\frac{S_{k-j}}{a_{k-j}}\in(-\frac{y}{a_{k-j}}+\frac{a}{g_ka_{k-j}}
-\frac{\gep}{g_{k-j}a_{k-j}}, -\frac{y}{a_{k-j}}+\frac{a}{g_k
a_{k-j}}+\frac{\gep}{g_{k-j}a_{k-j}} )\Big),
\end{multline}
where in the last inequality we used the monotonicity of $g$.
Again \eqref{StoneLocal} applies, and in conjunction with
\eqref{localEst} gives

\begin{multline*}\sup\limits_{y\in((a-\gep)/g_j, (a+\gep)/g_j)}\PPP\Big(
S_{k-j}\in(-y+\frac{a-\gep}{g_k}, -y+\frac{a+\gep}{g_k} )\Big)\\ <
cf_\ga (0) \frac{2\gep}{a_{k-j} g_{k-j}}<c^2\,
 \PPP(A_{k-j}(a, \gep))
 \end{multline*}
for  $k, j\ge n_0$ with $k-j \ge n_0$. Relation \eqref{localEst}
implies that for every $\gep>0$ the two series
$\sum\nolimits_{n=0}^{+\infty}\PPP(A_n(a, \gep))$ and
$\sum\nolimits_{n=1}^{+\infty}(a_n g_n)^{-1}$ converge or diverge
together. When $\sum\nolimits_{n=1}^{+\infty}(a_n
g_n)^{-1}<\infty$, the first Borel-Cantelli lemma applied to the
sequence of the $A_n(a, \gep)$ 's with the choices $a=0$,
$\gep=k$, $k\in \mathbb{N}$, gives $\lim_{n\to+\infty} g_n
|S_n|=+\infty$ (i.e., there are no finite accumulation points).
When this series diverges, the result follows on applying Lemma
\ref{BC} to the $A_n(a, \gep)$ 's with the choices $a\in\D{Q}$,
$\gep=1/k$, $k\in \mathbb{N}\setminus\{0\}$ and observing that
$K(gS)$ is a closed subset of $\ol{\D{R}}$.

\textbf{Step 2:} For the sequence $g_n:=n$, using Facts 1,2, we
see that \eqref{CramerControl} holds and
$\sum\nolimits_{n=1}^{+\infty}\frac{1}{a_n g_n}<+\infty$.
Consequently, by Step 1, $K(gS)=\{-\infty, +\infty\}$.

\textbf{Step 3:} Now let $(g_n)_{n\ge1}$ be any increasing
sequence of positive numbers. Introduce $g_n^\prime=\min\{g_n,
n\}$ for $n\ge1$. Then $(g_n^\prime)_{n\ge1}$ satisfies
$\lim_{n\to+\infty}g_n^\prime a_n(\theta+\gep)^n=0$.

If $\sum\nolimits_{n=1}^{+\infty}(a_n g_n)^{-1}<\infty$, then
since $1/g_n^\prime<1/g_n+1/n$ and
$\sum\nolimits_{n=1}^{+\infty}(a_n n)^{-1}<\infty$, it follows
that $\sum\nolimits_{n=1}^{+\infty}(a_n g_n^\prime)^{-1}<\infty$.
By the result of Step 1, it follows that $\varliminf_{n\to+\infty}
g_n^\prime |S_n|=+\infty$. Combining this with $g_n^\prime\le
g_n$, we get $\varliminf_{n\to+\infty} g_n|S_n|=+\infty$.

If $\sum\nolimits_{n=1}^{+\infty}(a_n g_n)^{-1}=+\infty$, then
$g_n^\prime\le g_n$ implies that
$\sum\nolimits_{n=1}^{+\infty}(a_n g_n^\prime)^{-1}=+\infty$. By
Step 1, $K(g^\prime S)=\ol{\D{R}}$; and since $\lim_{n\to+\infty}
n |S_n|=+\infty$ (Step 2), we conclude that $K(gS)=\ol{\D{R}}$.
  \hfill \qedsymbol

\medskip

\nid \textbf{Proof of Proposition \ref{propRec}:}

(a) We set $A_n(\gep):=[|Y_n|<\gep]$ for all
$n\in\D{N}\setminus\{0\}$ and $\gep>0$. By the assumption,
$\PPP(\cup_{k=n}^{+\infty} A_k(\gep))=1$ for all $\gep>0$ and
$n\ge1$. Recursively we construct a strictly increasing sequence
of positive integers $(n_j)_{j\ge0}$ with $n_0=1$ and
\begin{equation}\label{BCRec}
 \PPP(\cup_{k=n_{j-1}}^{n_j-1}
A_k(\frac{1}{j}))\ge 1-\frac{1}{2^j}.
\end{equation}
Define $g_n=\sqrt{j}$ for $n_{j-1}\le n < n_j$. By \eqref{BCRec}
and the first Borel-Cantelli lemma, we have that with probability
one there is a random $j_0$ so that for all $j\ge j_0$ there is a
$k_j\in [n_{j-1}, n_j-1]$ with $|Y_{k_j}|<1/j$; i.e.,
$g_{k_j}|Y_{k_j}|<1/\sqrt{j}$. This proves (a).

(b) For every $n\ge1$, there is a $\gd_n\in(0, 1)$ so that
$\PPP(|Y_n|<\gd_n)\le 2\max\{P(Y_n=0), 2^{-n}\}$. Let
$h_n:=n/\gd_n$ for all $n\ge1$. Then $\sum_{n=1}^{+\infty}
\PPP(h_n|Y_n|<n)<+\infty$, which implies that
$\lim_{n\to+\infty}h_n|Y_n|=+\infty$ a.s. Let $g_n=\max\{h_1,
h_2,\cdots, h_n\}$. \hfill \qedsymbol

\medskip

Proposition 1(a) says that all recurrent sequences $(Y_n)_{n\ge1}$
can afford to get increased and still be recurrent. There is no
recurrent sequence that is ``on the brink of loosing its
recurrence.''

\section{The multidimensional case}
 Again we consider a sequence
$(X_i)_{i\ge1}$ of i.i.d random variables with values in $\D{R}^d$
with $d\ge2$ and the corresponding random walk $(S_n)_{n\ge1}$
with $S_n=X_1+\cdots+X_n$. If the support of the distribution of
$X_1$ is not contained in a two dimensional subspace of $\D{R}^d$,
then the walk is transient (see \cite{DU}, Chapter 3, Theorem
2.12). So the only meaningful multidimensional case is when $d=2$.
Our results, concerning random walks in $\D{R}$, have analogs in
this case too. Now for a sequence $x:=(x_n)_{n\ge1}$ in $\D{R}^2$,
we define $$K(x):=\{a\in\D{R}^2\cup\{\infty\}:a \text{ is an
accumulation point for } (x_n)_{n\ge1}\}$$
 where $\D{R}^2\cup\{\infty\}$ is the one point compactification of
 $\D{R}^2$ equipped with the usual topology. In this topology, $K(x)$ is a closed set.

The analog of Theorem \ref{symmetric} holds in dimension two also.
We only need in its statement to write $\|S_n\|_{\infty}$ instead
of $|S_n|$ (of course any other equivalent norm in $\D{R}^2$
works). The crucial element in its proof is again the analog of
Lemma \ref{dominance}, which reads $\PPP(\|X-x\|_{\infty}\le
\gd)\le 16\ \PPP(\|X\|_{\infty}<\gd)$ for any two dimensional
random variable $X$ with nonnegative characteristic function, and
any $x\in\D{R}^2, \gd>0$. The proof is done in the same way as
that of Lemma \ref{dominance} using the function $A((x_1,
x_2)):=\int_{[0,1]\times[0,1]} \PPP(|X_1-x_1|<r_1, |X_2-x_2|<r_2)\
dr_1 dr_2$, where $X=(X_1, X_2)$.

For the analog of Theorem \ref{stable}, we assume that $F$, the
common distribution function of the $X_i$'s, has support not
contained on a proper subspace of $\D{R}^2$. Also that it belongs
to the domain of attraction of a stable distribution $F_\ga$ with
index $\ga\in(0,2]$, with zero centering constants, and that the
characteristic function $ \varphi$ of $F$ (defined by
$\varphi(t):=\int_{\D{R}}e^{it\cdot x}dF(x)$ for all
$t\in\D{R}^2$) satisfies Cramer's condition (C); that is,
$\varlimsup_{|t|\to+\infty}|\varphi(t)|<1$. Then it can be seen,
using the local limit theorem of Stone and well known criteria for
recurrence (see, e.g., \cite{DU}, Chapter 3, Lemma 2.4), that the
walk is recurrent only when $\ga=2$. The next theorem concerns
this case.

\begin{theorem} \label{2dstable}
Assume that $F$ belongs to the domain of attraction of a
non-degenerate two dimensional normal distribution, with zero
centering constants, and moreover that the characteristic function
of $F$ satisfies Cramer's condition (C). Then for any increasing
sequence $(g_n)_{n\ge 1}$ of positive numbers, it holds
\begin{equation*}
K(gS)=
    \begin{cases} \{\infty\} & \text{a.s. if } \sum\limits_{n=1}^{+\infty}\frac{1}{(a_n g_n)^2}<+\infty,  \medskip \\
\D{R}^2\cup\{\infty\}  & \text{a.s. if
}\sum\limits_{n=1}^{+\infty}\frac{1}{(a_n g_n)^2} =+\infty.
    \end{cases}
\end{equation*}
\end{theorem}

\nid The proof goes exactly as that of Theorem \ref{stable}.
\bigskip

\nid \textbf{Acknowledgement.} I am grateful to B\'alint Vir\'ag
for showing me a proof of Lemma \ref{dominance} before I got the
one included now. Also for comments that improved the presentation
of the paper.

\bibliographystyle{annals}

\bibliography{bibliography}

\begin{thebibliography}{8}
\providecommand{\natexlab}[1]{#1}
\providecommand{\url}[1]{{\tt #1}}
\providecommand{\urlprefix}{URL }
\expandafter\ifx\csname urlstyle\endcsname\relax
  \providecommand{\doi}[1]{doi:\discretionary{}{}{}#1}\else
  \providecommand{\doi}{doi:\discretionary{}{}{}\begingroup
  \urlstyle{rm}\Url}\fi
\providecommand{\bibAnnoteFile}[1]{%
  \IfFileExists{#1}{\begin{quotation}\noindent\textsc{Key:} #1\\
  \textsc{Annotation:}\ \input{#1}\end{quotation}}{}}
\providecommand{\bibAnnote}[2]{%
  \begin{quotation}\noindent\textsc{Key:} #1\\
  \textsc{Annotation:}\ #2\end{quotation}}

\bibitem[{Breiman(1992)}]{BREI}
Breiman, L. (1992). Probability, Society for Industrial and Applied
  Mathematics, Philadelphia, PA.
\bibAnnoteFile{BREI}

\bibitem[{Chung and Erdos(1947)}]{CE}
Chung, K.L. and Erdos, P. (1947). On the lower limit of sums of independent
  random variables, Ann. Math. {\bf 48}, no.~4: 1003--1013.
\bibAnnoteFile{CE}

\bibitem[{Durrett(1996)}]{DU}
Durrett, R. (1996). Probability: Theory and examples, Wadsworth Pub. Co, second
  edn.
\bibAnnoteFile{DU}

\bibitem[{Feller(1971)}]{F}
Feller, W. (1971). An introduction to probability theory and its applications,
  vol.~2, John Wiley and Sons, Inc., New York-London-Sydney, second edn.
\bibAnnoteFile{F}

\bibitem[{Ibragimov and Linnik(1971)}]{IB}
Ibragimov, I.~A. and Linnik, Yu.~V. (1971). Independent and stationary
  sequences of random variables. With a supplementary chapter by I. A.
  Ibragimov and V. V. Petrov. Translation from the Russian edited by J. F. C.
  Kingman, Wolters-Noordhoff Publishing, Groningen.
\bibAnnoteFile{IB}

\bibitem[{Petrov(1979)}]{PE}
Petrov, V.~V. (1979). Remark on the lower limit for the modulus of sums of
  independent random variables, Lithuanian Math. J. {\bf 18}, no.~4: 528--531.
\bibAnnoteFile{PE}

\bibitem[{Petrov(1995)}]{PE2}
Petrov, V.~V. (1995). Limit theorems of probability theory. Sequences of
  independent random variables, Oxford University Press, New York.
\bibAnnoteFile{PE2}

\bibitem[{Stone(1965)}]{ST}
Stone, C. (1965). A local limit theorem for non-lattice multidimensional
  distribution functions, Ann. Math. Statist {\bf 36}: 546--551.
\bibAnnoteFile{ST}

\end{thebibliography}

\end{document}